\documentclass[12pt]{article}
\usepackage{amsmath,amsthm,amsfonts}
\theoremstyle{plain}

\theoremstyle{plain}

\theoremstyle{remark}

\theoremstyle{remark}

\theoremstyle{plain}

\normalbaselines
\begin{document}
\renewcommand{\theequation}{\arabic{equation}}

\def\v{\forall}
\def\ex{\exists}
\def\w{W_{\alpha }}
\def\eps{\varepsilon}

\title{RADEMACHER CHAOS IN SYMMETRIC SPACES}
\author{\sc S.\,V. Astashkin}


\date{}
\maketitle

\begin{abstract}
In this paper we study properties of series with respect to orthogonal
systems $\{r_i(t)r_j(t)\}_{i\not = j}$ and $\{r_i(s)r_j(t)\}_{i,j=1}^\infty $
in symmetric spaces on interval and square, respectively.
Necessary and sufficient conditions for the equivalence of these systems
with the canonical base in $l_2$ and also for the complementability of the
corresponding generated subspaces, usually called {\it Rademacher chaos}, are
derived. The results obtained allow, in particular, to establish the
unimprovability of the
exponential integrability of functions from Rademacher chaos. Besides, it is
shown that for spaces that are "close"  to $L_\infty $, the systems considered,
in contrast to the ordinary Rademacher system, do not possess the property of
unconditionality. The degree of this non-unconditionality is explained in the
space $L_\infty $.
\end{abstract}

\section{Introduction}

\def\ab{I\times{I}}
\def\bc{\{r_i(s)r_j(t)\}_{i,j=1}^{\infty}}
\def\cd{\{r_i(t)r_j(t)\}_{i<j}}
\def\de{L_{\infty}}
\def\ef{\sum_{i,j=1}^{\infty}a_{i,j}r_i(s)r_j(t)}
\def\fg{\sum_{j=1}^{\infty}a_{i,j}r_j(t)}
\def\gh{\chi_{(0,t)}}
\def\hi{\ln(e/t)}
\def\ij{\ln(e/s)}
\def\jk{a_{i,j}(x)}
\def\kl{a_{i,j}(y)}
\def\lm{T_{u,v}}
\def\no{Q'_{i,j}}
\def\mn{Q'_{i,j}(T_{u,v}x)\,du\,dv}
\def\op{Q_{i,j}}
\def\pr{x(u,v)\,du\,dv}
\def\rs{T_{u,v}x\,du\,dv}
\def\ls{\int_{U_i}}
\def\lv{\int_{{\bar U}_i}}
\def\lg{\int_{U_j}}
\def\lx{\int_{{\bar U}_j}}
\def\st{U_i\times{U_j}}
\def\tv{\bar{U}_i\times{U_j}}
\def\vw{U_i\times{\bar{U}_j}}
\def\wz{\bar{U}_i\times{\bar{U}_j}}
\def\zx{\int_0^1\int_0^1x(u,v)r_{i+1}(u)r_{j+1}(v)\,du\,dv}
\def\xy{\theta_{i,j}}
\def\ya{\sum_{i,j=1}^{2^k}\theta_{i,j}r_i(s)r_j(t)}
\def\zb{\sum_{i=1}^{2^k}\theta_{i,j}r_i(s)}
\def\xc{\sum_{j=1}^{2^k}}
\def\oh{\sum_{i=1}^{2^k}}
\def\nx{M(\varphi_{\varepsilon})}
\def\wg{T_{\theta}}
\def\fv{\sum_{j=1}^{\infty}c_jr_j}
\def\ip{\Lambda_p(\varphi)}
\def\pq{(a_{i,j})_{i,j=1}^{\infty}}
\def\bs{\sum_{i\ne j}b_{i,j}r_i(t)r_j(t)}
\def\bt{\sum_{i\in D,j\not\in D}a_{i,j}}
\def\bj{\sum_{i\in D,j\in E}a_{i,j}}
\def\bi{\varphi_n(\theta)}
\def\bk{2^{-n^2}\,\sum_{\theta}\varphi_n(\theta)}
\def\bl{\biggl|\biggl|\sum_{i,j=1}^n{\mbox r}_{i,j}(u)r_i(s)r_j(t)\biggr|\biggr|_{\infty}}
\def\ck{\sum_{i=1}^n\theta_{i,j}r_i(s)}
\def\cl{\sum_{j=1}^n}

\def\nl{\vskip 0.3cm}
\def\be{\begin{multline*}}
\def\ee{\end{multline*}}
\def\ds{\;}

\nl Let
$$
r_k(t)=\,{\rm sign}\sin{2^{k-1}{\pi}t}\,\,(k=1,2,\dots )
$$
be a Rademacher system on $I=[0,1]$. The set of all real functions
$y(t)$ which can be presented in the form
$$
y(t)=\,\bs\,\quad (t\in I) \leqno(1)
$$
is called the {\it
Rademacher chaos } of degree 2 with respect to this system.
The orthonormalized system $\cd$ on $I$, unlike the ordinary
Rademacher system,
consists of functions which are not independent.
Nevertheless, its properties remind in many aspects the properties of a family of
independent and uniformly bounded functions.
This concerns, in particular, the integrability of the functions of form
(1): The condition $\sum_{i,j}b_{i,j}^2\,<\,\infty$
implies the summability of the function
$\exp(\alpha|y(t)|)$ for each  $\alpha>0$
\cite[p.105]{1}. In the same time, there are substantial differences.
For example, the system $\cd$ is not a Sidon system
(see \cite{2}).
\par
The main purpose of this paper is a study of the behaviour of the
Radema\-cher chaos in arbitrary symmetric spaces of functions defined on
an interval. This will allow us to specify essentially the above
mentioned results and to obtain new assertions about the geometric structure of
subspaces of the symmetric spaces consisting of functions of the form (1). The
essence of the method we use in the sequel is the passage to the so called
"decoupling" chaos, that is, the set of functions, defined on the square
$\ab$, of the form
$$
x(s,t)=\,\ef.\leqno{(2)}
$$

Note that the study of the decoupling Rademacher chaos
(i.e., multiple series with respect to this system)
is of interest by itself.  \vskip 0.2cm

Let us recall that a Banach space $X$
of Lebesgue measurable on $I$ functions
$x=x(t)$ is called symmetric (s.s.), if:

1) The condition $|x(t)|\le{|y(t)|},~y\in X$ implies  $x\in X$ and
$||x||\le{||y||};$

2) The assumption
$\mu\{t\in I:\,|x(t)|>\tau\}=\mu\{t\in I:\,|y(t)|>\tau\}\,(\tau> 0)$
($\mu$ is the Lebesgue measure on $I$), and  $y\in X$, implies:
$x\in X$ and $||x||=||y||$.

Given any s.s. $X$ on the interval $I$, one can construct a space
$X(\ab)$ on the square $\ab$, having properties which are similar to
1) and	2). Indeed, if $\Pi:\,I\to{\ab}$
is one-to-one mapping that preserves the measure,
then  $X(\ab)$ (we shall call it s.s. too)
consists of all
Lebesgue measurable on	$\ab$ functions $x=x(s,t)$ for which
$x(\Pi(u))\in X$ and  $||x||_{X(\ab)}\,=\,||x(\Pi)||_X$.
\vskip 0.2cm

In the first part of the paper we obtain necessary and sufficient conditions
for the equivalence of the systems
$\cd$ and  $\bc$ (in what follows we call them the multiple Rademacher
systems)
with the canonical base in $l_2$, and also, for the complementability of the
subspaces generated by them.
The theorems proved resemble, in the essence,
the analogous results about usual
Rademacher series
in s.s. on an interval (see \cite{3}, \cite{4}, \cite[2.b.4]{5}). In the same
time, there are some differences:  The role of the "extreme" space is played
by the space $H$, the closure of $\de$ in the Orlicz space
$L_M$, $M(t)=e^t-1$,
instead of $G$, the closure in the Orlicz space
$L_N$, $N(t)=e^{t^2}-1$ (see Theorems 1 -- 4).
\par
The behaviour of the multiple Rademacher systems "close" to the space
$\de$, and in particular,  in the same space, becomes more complicated.
It is well-known \cite{6} that the usual Rademacher system
is a symmetric basic sequence in every s.s.
 $X$, that is, for arbitrary numerical sequence
$(c_j)_{j=1}^{\infty}$, we have
$$
\biggl|\biggl|\fv\biggr|\biggr|_X\,=\,\biggl|\biggl|
\sum_{j=1}^{\infty}c_j^*r_j\biggr|\biggr|_X
$$
($(c_j^*)_{j=1}^{\infty}$ is the rearrangement
of $(|c_j|)_{j=1}^{\infty}$
in decreasing order). On the contrary, the multiple systems
in  $\de$ and in the "close" spaces do not possess even the weaker
unconditionality property. Another example:
many s.s. sequences, which are intermediate between  $l_1$
and  $l_2$, are spaces of coefficient sequences
of Rademacher series from functional s.s. on the interval
(\cite{3}, \cite{7}).
This is not true for multiple systems, for example
in the case of the space $l_p~(1\le p<2)$ (see Theorems  5 -- 8
and the corollaries from them). \vskip 0.2cm
\par
Let us recall certain definitions. If  $X$ is a s.s.
on
$I$, then we shall denote by  $X^0$ the closure of $\de$ in  $X$.
For any two measurable on $I$ functions $x(t)$ and
$y(t)$ we set:
$$
<x,y>\,=\,\int_0^1x(t)y(t)\,dt
$$
(if the integral exists). The last notation is similarly
interpreted in the case of functions, defined
in the square $\ab$.
The {\it associated } space  $X'$ with $X$ is defined as the space
of all measurable functions  $y(t)$ for which
the norm $$
||y||_{X'}\,=\,\sup\{<x,y>:\,||x||_X\le 1\}
$$
is finite.
It is not difficult to verify that for $x\in{\de}$,
$||x||_X=||x||_{X''}$, and thus $X^0=(X'')^0$
 \cite[p.255]{8}.

The norm in s.s. $X$ is said to be {\it order semi-continuous},
if the conditions $x_n=x_n(t)\ge 0(n=1,2,\ldots )$, $x_n\uparrow{x}$
almost everywhere on $I$, $x\in{X}$, imply: $||x_n||_X\to{||x||_X}$.
If the norm in	$X$ is order semi-continuous, then it is isometrically
embedded in  $X''$  \cite[p.255]{8}, that is,
$$
||x||_X\,=\,\sup\{<x,y>:\,||y||_{X'}\le 1\}.
$$
Important examples of s.s. are the Orlicz spaces
$L_S$ ($S(t)\ge 0$ is a convex and continuous function
on $[0,\infty)$) with the norm
$$
||x||_{L_S}\,=\,\inf\left\{u>0:\,\int_{0}^{1}\,S(|x(t)|/u)\,dt
\,\le{\,1}\right\}
$$
and the Marcinkiewicz space  $M(\varphi)$
($\varphi(t)\ge 0$ is a concave increasing function
on $[0,1])$ with the norm
$$
||x||_{M(\varphi)}\,=\,\sup\left\{\frac{1}{\varphi(t)}\,
\int_o^t\,x^*(s)\,ds:\,0<t\le 1\right\}
$$
($x^*(s)$ is the decreasing left continuous
rearrangement of the function  $|x(t)|$  \cite[p.83]{9}).

\section { Equivalence of the multiple Rademacher systems
to the canonical base in $l_2$}

\nl {\bf Theorem~1.} {\it The system $\bc$ in the s. s.  $X(\ab)$
is equivalent to the canonical base
in $l_2$, if and only if  $X\supset H$ where
$H=L_M^0$, $M(t)=e^t-1$.
\par
\medskip
Proof}. The equivalence of the system $\bc$
to the canonical base in  $l_2$ means that
for arbitrary numerical sequence
$a=\pq$, we have
$$
\bigg\|\ef\bigg\|_{X(\ab)}\,\asymp{\,||a||_2},\leqno{(3)}
$$
where  $||a||_2=(\sum_{i,j=1}^{\infty}a^2_{i,j})^{1/2}$
(that is, a two-sided estimate takes place with
constants that depend only on  $X$).
\par
Assume first that  $X\supset H$. For $a=\pq\in{l_2}$
denote
$$
x(s,t)=\,\ef.\leqno(4)
$$
If  $q\ge 2$, then after a second application of  Khinchin's inequality
 \cite[p.153]{10} and the Minkowski integral inequality
  \cite[p. 318]{11} we obtain:

$$
||x||_q\;=\;\left(\int_0^1\int_0^1|x(s,t)|^q\,ds \,dt \right)^{1/q}$$
$$
\le\;{\,\sqrt{q}\left\{\int_0^1\left(\sum_{i=1}^{\infty}|\fg|^2\right)^{q/2}\,dt\right\}^{1/q}}$$
$$
\le\;\sqrt{q}\left\{\sum_{i=1}^{\infty}\left(\int_0^1|\fg|^q\,dt\right)^{2/q}\right\}^{1/2}\,\le{\,q||a||_2}.
$$
The inequality	$||x||_q\le{q||a||_2}$ obviously holds also for
 $1\le q<2$.

Hence, making use of the expansion
$$
\exp(uz)\,=\,\sum_{k=0}^{\infty}\frac{u^k}{k!}z^k \quad (u>0),
$$
we estimate:
$$
\int_0^1\int_0^1\left\{\exp(u|x(s,t)|)-1\right\}\,ds\,dt$$
$$
=~\sum_{k=1}^{\infty}\frac {
u^k}{k!}\int_0^1\int_0^1|x(s,t)|^k\,ds\,dt
\le{\sum_{k=1}^{\infty}\frac{u^k}{k!}k^k||a||_2^k.}
$$
Now if $u<(2e||a||_2)^{-1}$, then the sum of the series
on the right hand side of the last inequality
does not exceed 1. Therefore, by the definition of the space
 $H$, the linear operator
$$
Ta(s,t)=\,\ef
$$
is bounded from $l_2$ in  $H(\ab)$, and, in addition,
$||Ta||_H=||x||_H\le{2e||a||_2}$. Since $X\supset H$, we derive
$||x||_X\le{2Ce||a||_2}$.

On the other hand, for the function $x$, defined by the relation
(4), in view of the embedding $X\subset{L_1}$ which holds for arbitrary
s.s. on  $I$  \cite[p.124]{9}, and also by virtue of  Khinchin's inequality
for the  $L_1$-norm with the constant from   \cite{12}, and Minkowski
inequality, we obtain :
$$
\begin{array}{rcl}
\|x\|_{X(\ab)}&\ge&D^{-1}\|x\|_1\\
&\ge&(\sqrt{2}D)^{-1}
\ds\int_0^1\biggl[\sum_{i=1}^\infty\biggl(\fg\biggr)^2\biggr]^{1/2}\,dt\\
&\ge&(\sqrt{2}D)^{-1}\biggl\{\ds\sum_{i=1}^\infty\biggl[\ds\int_0^1\biggl|\fg
\biggr | \,dt \biggr]^2\biggr\}^{1/2}\\
&\ge&(2D)^{-1}\|a\|_2
\end{array}
\leqno{(5)}
$$
Thus  (3) is proved.

For the proof of the inverse statement we need two lemmas.

In what follows
$$
n_x(z)=\,\mu\{\omega:\,|x(\omega)|>z\},\quad
x^*(t)=\,\inf\{z>0:\,n_x(z)<t\}\,
$$ are, respectively, the distribution function and the decreasing, left side
continuous rearrangement of the measurable on $I$ or on $\ab$ function
$|x(\omega)|$ \linebreak \cite[p.83]{9}.
\par
\medskip
{\bf Lemma 1.} {\it If	$n_{x_k}(z)\to{n_x(z)}$ for  $k\to{\infty}$,
then $x_k^*(t)\to{x^*(t)}$ at all points of continuity
of the function $x^*(t)$ (and hence, almost everywhere).
\par
\medskip
Proof}.  Let  $x^*(t)$ be continuous at the point  $t_0$
and  $\varepsilon>0$. For \linebreak $z_1=\,x^*(t_0)+\varepsilon$ we have
$n_x(z_1)<t_0$ and thus  $n_{x_k}(z_1)<t_0$,
if  $k\ge{k_1}$. Then, from the definition of the rearrangement,
$$
x_k^*(t_0)\le{x^*(t_0)+\varepsilon}\,\,(k\ge{k_1}). \leqno(6)
$$

On the other hand, by the continuity of $x^*(t)$ at  $t_0$,
there exists a	$\delta>0$ such that
$$
x^*(t_0+\delta)\ge{x^*(t_0)-\varepsilon}
\leqno{(7)}
$$
If  $z_2=\,x^*(t_0+\delta)-\varepsilon$, then
$n_x(z_2)\ge{t_0+\delta}$ and there exists
 $k_2$ such that   $n_{x_k}(z_2)\ge{t_0}$ for $k\ge{k_2}$.
Therefore, in view of (7),
$$
x_k^*(t_0)\ge{x^*(t_0)-2\varepsilon}\,\,(k\ge{k_2}).
$$
For  $k\ge{\max(k_1,k_2)}$, both the previous inequality and
inequality (6) hold.
Combining all them we complete the proof.
\vskip 0.3cm

Denote
$$
L(z)=\,\mu\{(s,t)\in{\ab}:\,\ln(e/s)\ln(e/t)>z\}\,(z>0).
$$
\par
\medskip
{\bf Lemma 2.} {\it For any  $z\ge 1$, we have
$$
1/2\exp(-2\sqrt{z}+2)\le{L(z)}\le{2\exp(-\sqrt{z}+2)}\leqno{(8)}
$$
\par
\medskip
Proof}.  If  $g_z=\,u+z/u,\,h_z=\,\max(z/u,u)$, then
$$
h_z(u)\,\le{\,g_z(u)}\,\le{\,2h_z(u)}\,(u>0). \leqno(9)
$$
After the change of variable $u=\,\ln(e/s)$ we obtain
$$
\begin{array}{rcl}
L(z)&=&\ds\int_0^1\mu\{t\in
I:\,\hi>z\ln^{-1}(e/s)\}\,ds\\ \\
&=&e\ds\int_0^1\exp\left(-\frac{z}{\ij}\right)\,ds
~=~e^2\ds\int_1^{\infty}\exp(-u-z/u)\,du\\ \\
&=&e^2\ds\int_1^{\infty}\exp(-g_z(u))\,du.
\end{array}
\leqno{(10)}
$$

Integrating by parts, we estimate
$$
\int_{\sqrt{z}}^ze^{-u}\frac{du}{u}=\frac{\exp(-\sqrt{z})}{\sqrt{z}}-\frac{\exp(-z)}{z}-\int_{\sqrt{z}}^ze^{-u}\frac{du}{u^2}\le{\frac{\exp(-\sqrt{z})}{\sqrt{z}}}.$$
Thus, in view of  (9) and  (10), changing the variable, we arrive at the
bound
$$
e^{-2}L(z)\;\le\;{\int_1^{\infty}\exp(-h_z(u))\,du}
~=~\int_1^{\sqrt{z}}\exp(-z/u)\,du
+ \exp(-\sqrt{z})$$
$$
=\;\exp(-\sqrt{z})+z\int_{\sqrt{z}}^ze^{-u}\,du/u^2$$
$$
\;\le\;{\exp(-\sqrt{z})+\sqrt{z}
\int_{\sqrt{z}}^ze^{-u}\,du/u}\le{2\exp(-\sqrt{z})}.
$$
For the reverse estimate, again with the help of
(9) and  (10), we obtain
$$
e^{-2}L(z)\ge{\int_1^{\infty}\exp(-2h_z(u))}\ge{\int_{\sqrt{z}}^{\infty}\exp(-2
u)\,du}=2^{-1}\exp(-2\sqrt{z}),
$$
and inequality (8) is proved.
\vskip 0.2cm

Returning now to the proof of Theorem 1, assume that the system \linebreak
$\bc$ is equivalent in $X(\ab)$ to the canonical base in
$l_2$. If  $a_{i,j}=b_ic_j$, then
$$
\bigg\|\ef\bigg\|_{X(\ab)}\,\le{\,C||(a_{i,j})||_2}\,=\,C||(b_i)||_2||(c_j)||_
2 .
$$
In particular, for  $b_i=c_i=1/\sqrt n\,(1\le i\le n)$	and
$b_i=c_i=0\,(i>n)$ \linebreak $(n=1,2,\ldots )$,
in view of the previous inequality and the symmetry of
$X$,
$$
||v_n^*(s)v_n^*(t)||_{X(\ab)}\,\le{\,C}\,(n=1,2,\ldots ), \leqno(11)
$$
where  $v_n(s)=\,1/\sqrt n\,\sum_{i=1}^nr_i(s)$.

By the central limit theorem
(see, for example,  \cite[p.161]{13}),
$$
\lim_{n\to{\infty}}\mu\{s\in I:\,v_n^*(s)>z\}\,=\,\Phi(z)
$$
where
$$
\Phi(z)=\,2/\sqrt{2\pi}\int_z^{\infty}\exp(-u^2/2)\,du\,(z>0).
$$
Applying Lemma 1, we then obtain
$$
\lim_{n\to{\infty}}v_n^*(s)\,=\,\Phi^{-1}(s)\,(s>0)
$$
($\Phi^{-1}(s)$ is the inverse function of $\Phi(z)$).
Therefore, in view of  (11) and the properties of the second associated space
$X''$  \cite[p.256]{8},
$$
||\Phi^{-1}(s)\Phi^{-1}(t)||_{X''(\ab)}\,\le C.
$$
Since $\Phi(z)\le{\exp(-z^2/2)}$, then
$\ln^{1/2}(e/s)\ln^{1/2}(e/t)\in{X''(\ab)}$ and thus, by Lemma 2,
in view of the symmetry of $X''$, we have
$$
g(t)=\,\hi\in{X''}.
$$
\par
If now	$\gh$ is the characteristic function
of the interval $(0,t)$, then, as known (see  \cite[p.88]{14}),
$$
h(t)=\,1/||\gh||_{L_M}=M^{-1}(1/t)=\ln(1+(e-1)/t)\asymp{g(t)}.$$
Therefore  $h\in{X''}$ as well. For arbitrary function  $x\in{\de}$ and
all  $t\in I$, in view of   \cite[p.144]{9}, we have
$$
x^*(t)\le{1/t\int_0^tx^*(s)\,ds}\le{t^{-1}||x||_{L_M}||\gh||_{(L_M)'}}=||x||_{L_M}h(t).$$
Therefore, $||x||_X=||x||_{X''}\le{K||x||_{L_M}}$ where
$K=\,||h||_{X''}$. This means that
$H=L_M^0\subset{X^0}\subset X$ and Theorem 1 is proved.
\vskip 0.3cm

We shall prove a similar result for the
"undecoupling" Rademacher chaos.
\par
\medskip
{\bf Theorem 2.} {\it The system $\cd$ in the s.s. $X$
on $I$ is equivalent to the canonical base
in $l_2$ if and only if $X\supset H$.
\par
\medskip
 Proof}.  The idea of the proof is a passage to the
"decoupling" chaos and then application of Theorem 1.

Let the function  $y(t)$ be of the form  (1) and
$$
y^N(t)=\,\sum_{1\le{i\ne j}\le N}b_{i,j}r_i(t)r_j(t) \leqno{(12)}
$$
for $N=1,2,\ldots $
The next representation follows by  combinatorics arguments.
$$
y^N(t)=\,2^{1-N}\sum_{D\subset{\{1,\ldots , N\}}}\bigg(\sum_{i\in D,~j\not\in
D}a_{i,j}r_i(t)r_j(t)\bigg ),\leqno{(13)}$$
where  $a_{i,j}=b_{i,j}+b_{j,i}$ and the summation is taken
over all  $D\subset{\{1,\ldots , N\}}$	\cite[p.108]{1}.
Clearly the functions
$$
{\bt}r_i(t)r_j(t)\,(t\in I)\,\,\,\,\mbox{and}\,\,\,\,{\bt}r_i(s)r_j(t)\,(s,t\in
I)$$ are equimeasurable, that is, they have equal distribution functions.
Hence, if s.s.
 $X$ belongs to $H$, then it follows from Theorem 1 and equality
(13) that
$$
||y^N||_X\le{\,2^{1-N}C\sum_{D\subset{\{1,\ldots,N\}}}\bigg(\sum_{i\in
D,~j\not\in D}a_{i,j}^2\bigg )^{1/2}}\le{\,2C||a||_2}.$$
Going to the limit as $N\to{\infty}$, we obtain
$||y||_X\le{\,2C||a||_2}$.

On the other hand, it is well-known (see, for example,	\cite[p.149]{15}), that
for a certain  $C_1>0$, $||a||_2\le{\,C_1||y||_1}$.
Therefore, in view of the embedding
$X\subset{L_1}$,
the reverse inequality also holds:
$||a||_2\le{\,C_1D||y||_X}$.

For the proof of the inverse proposition we need the following.
\par
\medskip
{\bf Lemma 3.} {\it Let
$$
x(s,t)=\,\sum_{i\in A,j\in B}a_{i,j}r_i(s)r_j(t)\,(s,t\in I),
\leqno(14)
$$
where $A,B$ are finite subsets of the set of
natural numbers ${\cal N}$. Then the distribution function
$n_x(z)$ will not change after the replacement	of the summation sets
$A$ and  $B$ by the sets $D$ and  $E$ respectively, where
$|D|=|A|$, $|E|=|B|$ ($|G|$ is the number of elements
of the set  $G$).
\par
\medskip
Proof}.  Indeed, by the Fubini theorem
$$
n_x(z)=\,\int_0^1\mu\biggl\{t\in I:\,\biggl|\sum_{j\in B}\biggl(\sum_{i\in A}a_{i,j}r_i(s)\biggr)r_j(t)\biggr|\,>\,z\biggr\}\,ds.$$
But it follows from the definition of the Rademacher function
(see also   \cite{6})
that for a fixed $s\in I$ the set measure under the integral sign
will not change after the substitution of $B$ by  $E$,
provided  $|E|=|B|$.
\vskip 0.2cm

Continuing now the proof of Theorem 2, assume that
the system of functions $\cd$ is equivalent to the canonical
base in  $l_2$.

For every function $x(s,t)$ of the form (14) there exists
an equimeasurable, in absolute value, function
$$
y(t)=\,{\bj}r_i(t)r_j(t)
$$
such that  $|D|=|A|$ and  $|E|=|B|$. It suffices to take
$D=A$ and   $E$ so that  $E\cap A=\emptyset$, $|E|=|B|$.  Then, by Lemma 3,
the absolute value of the function
$$
z(s,t)=\,{\bj}r_i(s)r_j(t)
$$
is equimeasurable to  $|x(s,t)|$, and in view of the independence of the
Radema\-cher functions and the choice of
$D$ and $E$,  to  $|y(t)|\,(t\in I)$.

It follows from the assumptions and the presented reasonings
that for a certain  $C>0$ and all  $N=1,2,\ldots $,
$$
\biggl|\biggl|\sum_{i,j=1}^{N}a_{i,j}r_i(s)r_j(t)\biggr|\biggr|_{X(\ab)}\,\le{\
C||a||_2}.
$$ After a passage to the limit as
$N\to{\infty}$ we obtain
$$
\biggl|\biggl|\ef\biggr|\biggr|_{X(\ab)}\,\le{\,C||a||_2}.
$$

Since the inverse inequality takes place always,
the system  $\bc$ is also equivalent to the canonical base
in $l_2$. Therefore, by Theorem 1,  $X\supset H$ and
Theorem 2 is proved.

{\bf Remark 1.} If the function  $y(t)$ can be presented in the
form  (1), then, as we mentioned in the introduction,
the condition  $\sum_{i,j}b_{i,j}^2\,<\infty$ implies the
summability
 of the function  $\exp(\alpha|y(t)|)$
for each  $\alpha>0$
 \cite[p.105]{1}.
Theorem~2 gives not only a new proof of this assertion, it shows
also  its unimprovability in the class of all s.s.

\section {Complementability of the Rademacher chaos in symmetric
space}

\nl Let $\tilde{M}(u)$ be the function which is complementary
to  $e^u-1$, that is,  $\tilde{M}(u)=
=\sup_{v>0}(uv-M(v))$. It can be easily shown that
 $\tilde{M}(u)\asymp{u\ln u}$ as $u\to{\infty}$.
As known (see	\cite[p.97]{14}), the space
 $H^*=H'=L_M'=L_{\tilde M}$ is the conjugate one to  $H$  and  $L_{\tilde
M}'=L_M$.
\par
\medskip
{\bf Theorem 3.} {\it Let  $X$ be a s.s. with order  semi-continuous
norm. The subspace  ${\cal R}(X)$ of all functions from  $X(\ab)$,
which admit a presentation of the form	$(2)$, is complemented in this
space if and only if
$H\subset X\subset{H'}$.
\par
\medskip
Proof}. If  $H\subset X\subset{H'}$, then by duality,
 $H\subset{X'}$. Hence, in view of Theorem 1, for arbitrary
$a=\pq\in{l_2}$, we have
$$
\biggl|\biggl|\ef\biggr|\biggr|_{X(\ab)}\,\le{\,C||a||_2}
\leqno{(15)}
$$
and
$$
\biggl|\biggl|\ef\biggl|\biggr|_{X'(\ab)}\,\le{\,C||a||_2}.\leqno{(16)}
$$

With any given	$x=x(s,t)\in{X(\ab)}$ we associate the system
$$
a=(\jk)_{i,j=1}^{\infty},\,\jk=\,\int_0^1\int_0^1x(s,t)r_i(s)r_j(t)\,ds\,dt
\leqno{(17)}
$$
We shall show that the orthogonal projector
$$
Px(s,t)=\,\sum_{i,j=1}^{\infty}\jk r_i(s)r_j(t)
\leqno{(18)}
$$
is bounded in  $X(\ab)$.

Indeed, in view of (16),
$$
||a||_2^2\;=\;\sum_{i,j=1}^{\infty}(\jk)^2
=\int_0^1\int_0^1x(s,t)Px(s,t)\,ds\,dt$$
$$
\;\le\;{||x||_X||Px||_{X'}}\le{C||a||_2||x||_X}.
$$
Then it follows from relations (15) that
$$
||Px||_X\,\le{\,C||a||_2}\,\le{\,C^2||x||_X}\,(x\in X).
$$

Since the image of $P$ coincides with ${\cal R}(X)$, then this
subspace is complemented.

For the proof of the inverse statement we shall need the following.
\par
\medskip
{\bf Lemma 4.} {\it If the projector $P$ defined by equality
$(18)$ is bounded in the s.s. $X(\ab)$, then  $H\subset
X\subset{H'}$.
\par
\medskip
Proof}. Denote by $||P||$ the norm of the projector $P$ in
the space  $X(\ab)$. Since
$$
<Py,x>\,=\,<y,Px>\,(x\in X,y\in{X'}),
\leqno{(19)}
$$
then $P$ is bounded also in  $X'(\ab)$ with a norm not exceeding
$||P||$. Let the sequences $(\jk),(\kl)$ be defined as in
(17). Then, in view of the order semicontinuity
of $X$,
the equality (19), and the pairwise orthogonality of the functions
$r_i(s)r_j(t)\,(i,j=1,2,\ldots )$, for arbitrary
$x\in{X(\ab)}$ we have
$$
||Px||_X\;=\;\sup\{<Px,y>:\,||y||_{X'}\le 1\}\;=$$
$$
=\;\sup\{<Px,Py>:\,||y||_{X'}\le 1\}\;=$$
$$
=\;\sup\biggl\{\biggl|\sum_{i,j}\jk\kl\biggr|:||y||_{X'}\le 1\biggr\}\;\le$$
$$
\le\;{\sup\{||(\jk)||_2||(\kl)||_2:\,||y||_{X'}\le 1\}}.
$$
Applying now (5) to the subspace  $X'$ and the function  $Py$,
we obtain:
$$
||(\kl)||_2\,\le{\,C||Py||_{X'}}\,\le{\,C\,||P||\,||y||_{X'}}.
$$
Therefore, $||Px||_X\le{C\,||P||\,||(\jk)||_2}$,
that is,
$$
\biggl|\biggl|\ef\biggr|\biggr|_{X(\ab)}\,\le{\,C\,||P||\,||(a_{i,j})||_2.}$$

Just in the same manner one can show that
$$
\biggl|\biggl|\ef\biggr|\biggr|_{X'(\ab)}\,\le{\,CA||(a_{i,j})||_2}.$$
Since the reverse  inequalities always hold
 (see (5) ), then by Theorem 1 we conclude that  $X\supset H$ and
$X'\supset H$. Using the fact that
 $X\subset{X''}$, by the duality we obtain
$H\subset X\subset{H'}$.
 \vskip 0.3cm

Next we continue the proof of Theorem 3.
Assume that the subspace
${\cal R}(X)$ is complemented in  $X(\ab)$. In view of Lemma 4, it
suffices to
show that the projector  $P$ is bounded in this
space.

Let
$s=\sum_{i=1}^{\infty}s_i2^{-i},\,u=\sum_{i=1}^{\infty}u_i2^{-i}\,(s_i,u_i=0,1)$
be the binary expansion of the numbers $s$ and	$u$ from  $I=[0,1]$.
Following   \cite[p.159]{16} (see also	 \cite{4} and  \cite{17}), we set
$$
s\dot{+}u=\,\sum_{i=1}^{\infty}|s_i-u_i|2^{-i}.
$$
This operation transforms the interval	$I$, as well as the
square	$\ab$, in a compact Abelian group.

For $u,v\in I$ we define on $\ab$ the shift transforms
$$
\psi_{u,v}(s,t)=\,(s\dot{+}u,t\dot{+}v).$$
Since they preserve the Lebesgue measure on
$\ab$, the operators
$$
\lm x(s,t)=x\left(\psi_{u,v}(s,t)\right)
$$
act in the s.s.  $X(\ab)$ isometrically.

Introduce the sets
$$
U_i=\,\biggl\{u\in
I:\,u=\sum_{k=1}^{\infty}u_k2^{-k},\,u_i=0\biggr\},\quad\bar{U}_i=\,I\setminus{
U _ i } \quad (i=1,2,\ldots).$$
Since  $r_1(s\dot{+}u)=r_1(s)$ and
$$
r_i(s\dot{+}u)\,=\,\left\{
\begin{array}{rl}
r_i(s), & \mbox{~~if~~} u\in{U_{i-1}}\\
-r_i(s),& \mbox{~~if~~} u\in{\bar{U}_{i-1}}
\end{array} \right.\,(i=2,3,\ldots)\leqno{(20)}$$
then the subspace  ${\cal R}(X)$ is invariant with respect to
the transforms	$\lm$. Therefore, by Rudin's theorem
 \cite[p.152]{18}, there exists
a bounded projector  $Q$ from $X(\ab)$ to this subspace
which commutes with all operators
$\lm$ \linebreak $(u,v\in I)$. We shall show that  $Q=P$.

Let us present $Q$ in the form
$$
Qx(s,t)=\,\sum_{i,j=1}^{\infty}\op(x)r_i(s)r_j(t),
\leqno(21)
$$
where $\op$ are linear functionals in  $X(\ab)$. Since
$Q$ is a projector, then
$$
\op(r_k(s)r_m(t))\,=\,\left\{
\begin{array}{ll}
1, & k=i,m=j\\
0, & \mbox{otherwise }.
\end{array} \right.
\leqno(22)
$$

>From equality (20) and the property $\lm Q=Q\lm$,
we obtain
for  \linebreak $\no=Q_{i+1,j+1}$ and  $i,j=1,2,\ldots $
$$
\no(\lm x)\,=\,\left\{
\begin{array}{rl}
\no(x), & u\in{U_i},v\in{U_j} \mbox{~~or~~}
u\in{\bar{U}_i},v\in{\bar{U}_j}\\
-\no(x), & u\in{U_i},v\in{\bar{U}_j} \mbox{~~or~~}
u\in{\bar{U}_i},v\in{U_j}.
\end{array} \right.$$
Hence, taking into account that $\mu(U_i)=1/2\,(i=1,2,\ldots )$, we conclude
that
$$
\ls\lg\mn=\lv\lx\mn=\frac{1}{4}\no(x),\leqno{(23)}
$$
$$
\ls\lx\mn=\lv\lg\mn=-\frac{1}{4}\no(x).\leqno{(23')}
$$

The functionals $\op$ are bounded in  $X(\ab)$. Indeed, since
\linebreak $\bc$ is an orthonormal system in  $L_2(\ab)$, then by the virtue of
Khinchin's inequality
with the constant from	\cite{12}, and the embedding $X\subset{L_1}$,
$$
|\op(x)|\,\le{\,||Qx||_2}\,\le{\,2||Qx||_1}\,
\le{\,2D||Qx||_{X(\ab)}}\le{\,2D||Q||||x||_{X(\ab)}},
$$
where  $||Q||$ is the norm of the projector $Q$ in the space
$X(\ab)$. Therefore, in relations  (23) and
(23'),  the functional can be taken out from the integral sign.
\par
It is easily verified that for any $s\in I$,
$$
\{z\in I:\,z=s\dot{+}u,u\in{U_i}\}\,=\,\left\{
\begin{array}{ll}
U_i, & \mbox{~~if~~} s\in{U_i}\\
\bar{U}_i, & \mbox{~~if~~} s\in{\bar{U}_i},
\end{array} \right. $$
$$
\{z\in I:\,z=s\dot{+}u,u\in{\bar{U}_i}\}\,=\,\left\{
\begin{array}{ll}
\bar{U}_i, & \mbox{~~if~~} s\in{U_i}\\
U_i, & \mbox{~~if~~} s\in{\bar{U}_i}.
\end{array} \right.$$
Then, introducing the notations
$$
k_{i,j}=\ls\lg\pr,l_{i,j}=\lv\lg\pr,$$
$$
m_{i,j}=\ls\lx\pr,n_{i,j}=\lv\lx\pr,$$
$$
K_{i,j}(s,t)=\,\chi_{\st}(s,t),\,L_{i,j}(s,t)=\,\chi_{\tv}(s,t),$$
$$
M_{i,j}(s,t)=\,\chi_{\vw}(s,t),\,N_{i,j}(s,t)=\,\chi_{\wz}(s,t),$$
we can write
$$
\ls\lg\rs=k_{i,j}K_{i,j}+l_{i,j}L_{i,j}+m_{i,j}M_{i,j}+n_{i,j}N_{i,j},$$
$$
\lv\lg\rs=l_{i,j}K_{i,j}+k_{i,j}L_{i,j}+n_{i,j}M_{i,j}+m_{I.j}N_{i,j},$$
$$
\ls\lx\rs=m_{i,j}K_{i,j}+n_{i,j}L_{i,j}+k_{i,j}M_{i,j}+l_{i,j}N_{i,j},$$
$$
\lv\lx\rs=n_{i,j}K_{i,j}+m_{i,j}L_{i,j}+l_{i,j}M_{i,j}+k_{i,j}N_{i,j}.$$
Thus, the relations (23) and (23'), the linearity of $\no$,
and the equalities
$$
r_{i+1}(s)=\chi_{U_i}(s)-\chi_{\bar{U}_i}(s) \quad (i=1,2,\ldots )
$$
yield
\smallskip

\centerline{
$\begin{array}{l}
\no(x) \hbox{\rule{0cm}{0.8cm}} \\
~~~~=~\ds{\no\left((K_{i,j}+N_{i,j}-L_{i,j}-M_{i,j})\zx\right)}
\hbox{\rule{0cm}{0.8cm}} \\
~~~~=~\ds{\no\left(\zx\,r_{i+1}(s)r_{j+1}(t)\right)} \hbox{
\rule{0cm}{0.8cm}} \\
~~~~=~\ds{\zx.} \hbox{\rule{0cm}{0.8cm}}
\end{array}$
}
\medskip

Therefore
$$
\op(x)\,=\,\jk\quad(i,j=2,3,\ldots ).
$$
A similar reasoning  shows that the last equality
remains true also for \linebreak $i,j=1,2,\ldots $  Thus, in view of
the relations  (18) and  (21), we obtain  $Q=P$ and hence the theorem is
proved.
\par
\medskip
{\bf Theorem 4.} {\it Let  $X$ be a s.s. with an order
semi-continuous norm in  $I$.
The subspace ${\bar{\cal R}}(X)$ of all functions from	$X$ which can be
presented in the form
$(1)$, is complemented in  $X$ if and only if
$H\subset X\subset{H'}$.  }
\par
\medskip
The proof is similar to that in Theorem
3. Omitting the details we only remark that
instead of the projector  $P$ one should consider also the orthogonal
projector
$$
Sx(t)=\,\sum_{1\le i<j<\infty}\,\int_0^1x(u)r_i(u)r_j(u)\,du\,r_i(t)r_j(t),
$$
and instead of the transforms $\psi_{u,v}$ of the square
$\ab$, the transforms  \linebreak $\psi_u(s)=s\dot{+}u\,(u\in I)$
of the interval  $I$.

\section { Rademacher chaos in	$\de$ and in "close" s.s.}

\nl For any given $n\in{\cal N}$ and
$\theta=\{\xy\}_{i,j=1}^n,~\xy=\pm 1$, we introduce the quantity
$$
\varphi_n(\theta)=\,\biggl|\biggl|\sum_{i,j=1}^n\xy r_i(s)r_j(t)\biggr|\biggr|_{\infty}$$
($||\cdot||_{\infty}$ is the norm in the space
$\de(\ab)$).

It follows from the definition of the Rademacher functions
that  $\sup_{\theta}\varphi_n(\theta)$
is attained for  $\xy=1~(i,j=1,\ldots ,n)$ and it equals  $n^2$.
\par
\medskip
{\bf Theorem 5.} {\it  We have
$$
\inf_{\theta}\bi\,\asymp{\,\bk}\,\asymp{\,n^{3/2}}
\leqno{(24)}
$$
with constants which do not depend on  $n\in{\cal N}$.
\par
\medskip
Proof}. Note first that by virtue of  Khinchin's inequality
with the constant from	\cite{12}, for any
$\xy=\pm 1$, we have
$$
\biggl|\biggl|\sum_{i,j=1}^n\xy r_i(s)r_j(t)\biggr|\biggr|_{\infty}
\;=\;\sup_{0<s\le 1}\cl\biggl|\ck\biggr|\;\ge$$
$$
\;\ge\;{\cl\int_0^1\biggl|\ck\biggr|\,ds}\;\ge$$
$$
\;\ge\;{\frac{1}{\sqrt{2}}\cl\biggl\{\int_0^1\biggl(\ck\biggr)^2\,ds\biggr\}^{1/2
} } =1/\sqrt{2}n^{3/2}.
$$
Therefore, $\inf_{\theta}\bi\,\ge{\,1/\sqrt 2n^{3/2}}$.

To prove the opposite inequality,
note first that
$$
\bk\,=\,\int_0^1\bl \,du,
\leqno{(25)}
$$
where  $\{{\mbox r}_{i,j}\}_{i,j=1}^n$ are the first  $n^2$
Rademacher functions numbered in an arbitrary order.

Apply the known theorem about the distribution of the
$\de$-norm of polynomials
with random coefficients (see	\cite[p. 97]{19}) to the linear
space $B$ of functions of the form
$$
f(s,t)=\,\sum_{i,j=1}^na_{i,j}r_i(s)r_j(t)
$$
defined on  $\ab$. Since for every function $f\in B$ we have
$|f(s,t)|=||f||_{\infty}$
on the square $K\subset{\ab}$ with a measure  $\mu(K)\ge{2^{-2(n-1)}}$,
then for  $z\ge 2$,
$$
\mu\biggl\{u\in I:\,\bl\,\ge{\,3\biggl[\sum_{i,j=1}^n\ln(2^{2n-1}z)\biggr]^{1/2}\biggr\}}\,\le{\,2/z}.$$
Now, after some not complicated transformations for the functions
$$
X(u)=\,\bl$$
we arrive at the estimate
$$
\mu\{u\in I: X(u)\ge{3\sqrt{2}n^{3/2}\tau}\}\,\le{\,e^{-\tau^2}}\,(\tau\ge 2).$$
Therefore, in view of (25),
$$
\bk\;=\;||X||_1\,=\,3\sqrt{2}n^{3/2}\,\int_0^{\infty}\mu\{u\in
I:\,X(u)\ge{3\sqrt{2}n^{3/2}\tau}\}d\tau$$
$$
\;\le\;{\,3\sqrt{2}n^{3/2}\left(2+\int_2^{\infty}e^{-\tau^2}d\tau\right)}\,\le{\,
9 \sqrt{2}n^{3/2},}
$$
and the theorem is proved.
\par

The "probability" proof of Theorem 5 does not yield the specific
arrangement of the signs for which the exact lower bound in (24) is attained.
That is why we give additionally one more proposition where, for
simplicity, only the case
$n=2^k$ is considered.
\par
\medskip
{\bf Proposition.} {\it For any  $k=0,1,2,\ldots $ there exists
an arrangement of the signs  $\theta=\{\xy\}_{i,j=1}^{2^k}$
for which
$$
\varphi_{2^k}(\theta)\,\le{\,2^{3k/2}}.
$$
\par
\medskip
Proof}. Consider the Walsh matrix, that is, the matrix
constructed for the first  $2^k$ Walsh functions
$$w_1(t),w_2(t),\ldots ,w_{2^k}(t)$$
($w_1(t)=1,w_{2^i+j}(t)=r_{i+2}(t)w_j(t),$ $i=0,1,\ldots ;j=1,\ldots,  2^i$)
 \cite[p.158]{16}.
Since  on the intervals
$\Delta_i^k=((i-1)2^{-k},i2^{-k})$ $(1\le i\le{2^k})$
these functions are constant and equal	 $+1$ or $-1$,	then one can
determine the signs
$$
\xy=\,{\rm sign}\,{w_j(t)},\quad t\in{\Delta_i^k}.
$$

For each $s\in I$ consider the function
$$
x_s(u)=\,\oh r_i(s)\chi_{\Delta_i^k}(u)\,(u\in I).
$$
Since  $|x_s(u)|=1$, then
$||x_s||_2=1$. For  $1\le j\le{2^k}$ the Fourier -- Walsh coefficients
of the function  $x_s(u)$  are given by
$$
c_j(x_s)=\,\oh\int_{\Delta_i^k}w_j(u)\,du\,r_i(s)\,=\,2^{-k}\oh \xy r_i(s).
$$
Then, by virtue of H\"older's and Bessel's inequalities, for the chosen
$\xy$ we have
$$
\biggl|\biggl|\ya\biggr|\biggr|_{\infty}\;=\;\sup_{0<s\le
1}\xc\biggl|\zb\biggr|=2^k\sup_{0<s\le 1}\xc|c_j(x_s)|\;\le$$
$$
\;\le\;{2^{3k/2}\sup_{0<s\le 1}\biggl\{\xc(c_j(x_s))^2\biggr\}^{1/2}}
\;\le\;{2^{3k/2}||x_s||_2}\,=\,{2^{3k/2}},
$$
and the proposition is proved.
\vspace{0.4cm}

For the
"non-decoupling" chaos, introduce the following quantity which is analogous
to $\bi$.
$$
{\bar\varphi}_n(\theta)=\,\biggl|\biggl|\sum_{i,j=1}^n\xy r_i(t)r_j(t)\biggr|\biggr|_{\infty}\leqno{(26)}$$
($\theta=\{\xy\}_{i,j=1}^n$ is the symmetric arrangement of the signs,
i.e.,  $\xy=\theta_{j,i}$, $||\cdot||_{\infty}$ is the norm in
$\de$ on  $I$).
\par
\medskip
{\bf Theorem 6.} {\it With a certain  $C>0$, independent of
$n\in{\cal N}$,
$$
\inf_{\theta}{\bar\varphi}_n(\theta)\,\asymp{\,2^{-n^2}
\sum_{\theta}{\bar\varphi}_n(\theta)}\,\asymp{\,n^{3/2}}
\leqno{(27)}
$$
\par
\medskip
Proof}. Observe that for arbitrary
$n\in{\cal N}$ and  $a=(a_{i,j})_{i,j=1}^n$
the following inequality holds:
$$
\biggl|\biggl|\sum_{i,j=1}^na_{i,j}r_i(t)r_j(t)\biggr|\biggr|_{\de}\,
\le{\,\biggl|\biggl|\sum_{i,j=1}^na_{i,j}r_i(s)r_j(t)\biggr|\biggr|_{\de(\ab)}},
$$
which, by Theorem 5, yields the inequalities  $\le$ in relation  (27).

Let $\theta$ be an arbitrary symmetric arrangement of the signs.
By Theorem~5, there exists a constant $C>0$ such that
$\|x\|\ge Cn^{3/2}$
for all natural $n$ and the functions
$$
x(s,t)=\,\biggl|\biggl|\sum_{i,j=1}^n\xy r_i(s)r_j(t)\biggr|\biggr|_{\infty}.
$$
Denote
$$
x_1(s,t)=\,\sum_{1\le i<j\le n}\xy r_i(s)r_j(t),\,x_2(s,t)=\,\sum_{1\le j<i\le n}\xy r_i(s)r_j(t),$$
$$
x_3(s,t)=x(s,t)-x_1(s,t)-x_2(s,t).
$$
Since  $\|x_3\|_{\infty}\le n$ and  $x_1(s,t)=x_2(t,s)$, then
$$
\|x_1\|_{\infty}\,\ge\,\frac{1}{4}Cn^{3/2}\leqno{(28)}$$
for all sufficiently large $n$.
Further we need a theorem from	\cite{20} about the comparison of the
distribution functions of quadratic and bilinear forms.
Let us formulate it: {\it Let
$X=(X_1,X_2,\ldots , X_n)$ be a vector with coordinates that are independent
symmetrically distributed random variables on a certain probability
space with a measure
$P$, and let
$Y_1,Y_2,\ldots , Y_n$ be independent copies of
$X_1,X_2,\ldots , X_n$, respectively.
Then there exist constants
$K_1,K_2,k_1,k_2$ such that for arbitrary forms
$$
Q(X)=\,\sum_{1\le i<j\le n}a_{i,j}X_iX_j,\,\,\tilde{Q}(X,Y)=\,\sum_{1\le i<j\le n}a_{i,j}X_iY_j$$
and any  $z>0$,
$$
K_1P\{k_1|Q(X)|>z\}\,\le\,P\{|\tilde{Q}(X,Y)|>z\}\,\le\,K_2P\{k_2|Q(X)|>z\}.
$$
}
Let us apply this theorem in the case
$X_i=r_i(t),Y_i=r_i(s)\,(t,s\in
I)$. Then, in view of inequality  (28) and the fact that
the norm of a function in  $L_{\infty}$ is defined by its distribution
function, we obtain
$$
\biggl|\biggl|\sum_{1\le i<j\le n}\xy r_i(t)r_j(t)\biggr|\biggr|_{\infty}\,
\ge C_1n^{3/2}
$$
with certain $C_1>0$ and for all sufficiently large
$n$. Since the arrangement of the signs is symmetric and the norm of the
diagonal terms in the sum do not exceed
$n$, we conclude that
$$
\biggl|\biggl|\sum_{i,j=1}^n\xy r_i(t)r_j(t)\biggr|\biggr|_{\infty}\,
\ge C_2n^{3/2}.
$$
Diminishing the constant $C_2$ one can make the last inequality hold
for all natural
$n$. The theorem is proved.
\par
\medskip
{\bf Remark 2.} The relation (27), obviously remains true also in the case
when the summation in
(26) is expanded only over $i<j$.
\vskip 0.3cm

Recall that an orthonormal system of functions
$\{u_k(t)\}_{k=1}^{\infty}$, defined on a certain probability space,
is called a {\it Sidon system}
(see, for example,  \linebreak \cite[p.327]{16}), if for every generalized
polynomial   ${\cal
P}(t)=\,\sum_{k=1}^na_ku_k(t)$ with respect to this system holds the estimate
$$
C^{-1}\sum_{k=1}^n|a_k|\,\le{\,||{\cal P}||_{\infty}}\,\le{\,C\sum_{k=1}^n|a_k|,}$$
where the constant  $C>0$ does not depend on the polynomial
${\cal P}(t)$.
\par
\medskip
{\bf Corollary.} {\it The multiple systems  $\bc$ and  $\cd$ are not
Sidon systems on  $\ab$ and $I$, respectively.}
\par
\medskip
{\bf Remark 3.} The assertion in the last corollary concerning the
system
$\cd$ was proved in   \cite{2}.
\vskip 0.4cm

The basic sequence  $\{x_n\}_{n=1}^{\infty}$
in a  Banach space $X$ is said to be unconditional
if the convergence of the series
$\sum_{n=1}^{\infty}a_nx_n\,(a_n\in{\cal
R})$ in  $X$ implies the convergence in the same space of the series
$\sum_{n=1}^{\infty}\theta_na_nx_n$
for arbitrary signs  $\theta_n=\pm 1\,(n=1,2,\ldots )$
 \cite[p.22]{16}. It is not difficult to verify that in this case
there exists a constant
$C>0$ such that
 $$
\bigg\|\sum_{n\in
F}a_nx_n\bigg\|_X\,\le{\,C\bigg\|\sum_{n=1}^{\infty}a_nx_n\bigg\|_X}
\leqno{(29)}
$$
for arbitrary  $F\subset{\cal N}$.
Besides, the smallest  $C$ for which (29) is true, is called
the {\it constant of unconditionality }
for the sequence  $\{x_n\}$.

As seen from Theorems 1 and 2, the basic sequences
$\bc$ and  $\cd$ are not only unconditional, they are symmetric
as well (see Introduction) in the s.s. $X(\ab)$ and  $X$, respectively,
provided  $X\supset H$. The situation is completely different in
the space $\de$ and in s.s. which are "close" to it.

We prove first one more auxiliary proposition.
\par
\medskip
{\bf Lemma 5.}{\it Let	$n_0=0<n_1<n_2<\cdots $; $c_{i,j}\in{\cal R}$ and
$$
y_k=\,\sum_{i,j=n_k+1}^{n_{k+1}}c_{i,j}r_i(s)r_j(t)\,\neq{\,0}.$$

Then  $\{y_k\}_{k=1}^{\infty}$ is an unconditional basic
sequence in any s.s. $X(\ab)$
with a constant of unconditionality equal to
$1$.
\par
\medskip
Proof}. It suffices to verify (see  \cite[p.23]{16}) that for arbitrary
$m=1,2,\ldots $, $\theta_k=\pm 1(k=1,2,\ldots ,m)$ and
real $a_k$, the norms of the functions \linebreak
$y=\,\sum_{k=1}^ma_ky_k$ and
$y_{\theta}=\,\sum_{k=1}^m\theta_ka_ky_k$ in  $X(\ab)$ coincide.

The values of the Rademacher functions
$r_i(s)\,(i=1,2,\ldots ,n_{m+1})$
give all possible arrangements of the signs
(up to multiplication by $-1$), corresponding to the intervals
$((k-1)2^{-n_{m+1}+1},k2^{-n_{m+1}+1})\,(k=1,2,\ldots ,2^{n_{m+1}-1})$.
That is why the distribution function
$$
n_{y_{\theta}}(z)\;=\;\mu\{(s,t)\in{\ab}:\,
|y_{\theta}(s,t)|>z\}\,=\,\int_0^1\mu\{s\in I:\,|y_{\theta}(s,t)|>z\}\,dt$$
$$
=\;\int_0^1\mu\biggl\{s\in
I:\,\biggl|\sum_{k=1}^m\theta_ka_k\sum_{i=n_k+1}^{n_{k+1}}\biggl(\sum_{j=n_k+1}^{n_{k+1}}c_{i,j}r_j(t)\biggr)r_i(s)\biggr|>z\biggr\}\,dt
$$
does not depend on the signs  $\theta_k$. Therefore
$||y_{\theta}||_{X(\ab)}=||y||_{X(\ab)}$ and the lemma is proved.
\vskip 0.3cm

For any arrangement of the signs
$\theta=\{\xy\}_{i,j=1}^{\infty},\xy=\pm
1$, define the operator
$$
\wg x(s,t)=\,\sum_{i,j=1}^{\infty}\xy a_{i,j}r_i(s)r_j(t)
$$
on  ${\cal R}(\de)$,
where
$$
x(s,t)=\ef\,\in{\,\de(\ab)},\,a=\pq\in{l_2}.
$$
\par
\medskip
{\bf Theorem 7.} {\it For arbitrary  $0<\varepsilon<1/2$,
there exists an arrangement of the signs
$\theta=\{\xy\}$, for which
$$
\wg:\,{\cal R}(\de)\,\not\to\,{\nx(\ab)},$$
where
$\nx$ is the Marcinkiewicz space defined by  the  function
$$\varphi_{\varepsilon}(t)=t\log_2^{-\varepsilon+1/2}(2/t)$$
{\rm(see Introduction)}.
\par
\medskip
Proof}.  By Theorem 5 and Lemma 3, for every
$k=0,1,2,\ldots $ one can find	 $\xy=\pm
1\,(2^k+1\le{i,j}\le{2^{k+1}})$ so that the associated function
$$
z_k(s,t)=\,\sum_{i,j=2^k+1}^{2^{k+1}}\xy r_i(s)r_j(t)$$
satisfies
$$
||z_k||_{\infty}\asymp{2^{3k/2}}.\leqno{(30)}$$
Set  $x_k(s,t)=\,2^{-(3+\varepsilon)k/2}z_k(s,t)$ and
$$
x(s,t)=\,\sum_{k=0}^{\infty}x_k(s,t)=\,\sum_{i,j=1}^{\infty}a_{i,j}r_i(s)r_j(t),$$
where  $a_{i,j}=\,2^{-(3+\varepsilon)k/2}\xy~(2^k+1\le{i,j}\le{2^{k+1}})$
and  $a_{i,j}=0$ otherwise.

It follows from (30) that
$$
||x||_{\infty}\,\le{\,C\sum_{k=0}^{\infty}2^{-\varepsilon k/2}}\,=
\,C2^{\varepsilon/2}/(2^{\varepsilon/2}-1),
$$
i.e.,  $x\in{{\cal R}(\de)}$.

Let the arrangement of the signs $\theta$ consists of the
values $\xy$\\$(2^k+1\le{i,j}\le{2^{k+1}})$ found above and arbitrary
$\xy=\pm 1$ in the other cases.
Then
$$
y(s,t)=\,\wg x(s,t)=\,\sum_{k=0}^{\infty}2^{-(3+\varepsilon)k/2}y_k(s,t)$$
where
$$
y_k(s,t)=\,\sum_{i,j=2^k+1}^{2^{k+1}}r_i(s)r_j(t).$$
Applying Lemma 5 to this sequence and
$X=\nx$, in view of the relation (29), we obtain
$$
||y||_{\nx}\,\ge{\,2^{-(3+\varepsilon)k/2}||y_k||_{\nx}}\quad(k=0,1,2,\ldots
).\leqno(31)
$$
It is clear from the definition of the Rademacher functions
that  $y_k(s,t)=2^{2k}$ for $0<s,t<
<2^{-2^{k+1}+1}$.
Thus, if  $u_k=\,2^{-2^{k+2}+1}\,(k=0,1,\ldots )$,
then the rearrangement satisfies  $y_k^*(u_k)\ge{2^{2k}}$.

Since, according to  \cite[p.156]{9},
$$
||x||_{\nx}\,\asymp{\,\sup_{0<u\le
1}x^*(u)\log_2^{\varepsilon-1/2}(2/u)},
$$
we obtain from the last inequality that
$$
||y_k||_{\nx}\,\ge{\,C2^{2k}\log_2^{\varepsilon-1/2}(2/{u_k})}\,\ge{\,C2^{(\varepsilon+3/2)k-1}}.$$
Then it follows from (31) that for every
$k=0,1,\ldots $,  $$||y||_{\nx}\,\ge{\,2^{\varepsilon k/2-1}}$$
and thus  $y=\,\wg x\,\not\in{\,\nx(\ab)}$.
\par
\medskip
{\bf Corollary 2.} {\it If $p\in{[1,2)}$,
then there is no a s.s.  $X$ for which
$$
\biggl|\biggl|\ef\biggr|\biggr|_{X(\ab)}\,\asymp{\,\biggl(
\sum_{i,j=1}^{\infty}|a_{i,j}|^p\biggr)^{1/p}}.\leqno{(32)}
$$
\par
\medskip
Proof}. If  $p=1$, then the assertion follows immediately from
Theorem 5, since for every s.s.
$X$ on	$I$,  $X\supset{\de}$ (see  \cite[p.124]{9}).

Assume that  (32) is true  for $p\in{(1,2)}$.
Then, taking  $a_{1,j}=c_j$ and $a_{i,j}=0~ (i\neq 1)$, we obtain
$$
\biggl|\biggl|\fv\biggr|\biggr|_X\,\asymp{\,\biggl(\sum_{j=1}^{\infty}|c_j|^p\biggr)^{1/p}}.$$
Therefore, by Theorem 3 from   \cite{3}, $X=\ip$
where  $\ip$ is the Lorentz space with the norm
$$
||x||\,=\,\biggl\{\int_0^1(x^*(t))^pd\varphi(t)\biggr\}^{1/p},\quad
\varphi(t)=\log_2^{1-p}(2/t).
$$

Since the fundamental function of this space
satisfies
$$||\gh||_{\ip}=\varphi^{1/p}(t)=\log_2^{-1+1/p}(2/t),$$
and the Marcinkiewicz
space is maximal among the s.s. with the same fundamental function
 \cite[p.162]{9},
then  $X\subset{\nx}$ where   $\varepsilon=1/p-1/2$.
Applying Theorem 7 we obtain  $\theta=\{\xy\},\xy=\pm
1$ and therefore
$$\wg:\,{\cal R}(\de)\not\to{\nx(\ab)}.$$
Moreover,
$$\wg:\,{\cal R}(\de)\not\to{X(\ab)},$$
and hence there exists a function  $x\in{{\cal
R}(\de)}\subset{X(\ab)}$ such that \linebreak $\wg{x}\not\in{X(\ab)}$.
Since this is in a contradiction with  the relation  (32),
the corollary is proved.  \vskip 0.3cm

Similar propositions hold in the
"undecoupling" chaos.

Recall that
${\bar{\cal R}}(L_{\infty})$ is a subspace of  $L_{\infty}$,
consisting of all functions of the form
$$
y(t)=\sum_{1\le i\ne j<\infty}b_{i,j}r_i(t)r_j(t),\,\,b=(b_{i,j})_{i,j=1}^{\infty}\in{l_2}.$$

For any arrangement of the signs
$\theta=\{\xy\}_{i,j=1}^{\infty},\xy=\pm
1$, we	define the operator
$$
\bar{T}_\theta y(t)=\,\sum_{i,j=1}^{\infty}\xy b_{i,j}r_i(t)r_j(t)
$$
in the space  ${\bar{\cal R}}(L_{\infty})$.
\par
\medskip
{\bf Theorem 8.} {\it For arbitrary  $0<\varepsilon<1/2$
there exists arrangement of the signs
$\theta=\{\xy\}$ for which
$$
\bar{T}_\theta :\,{\bar{\cal R}}(\de)\not\to{\nx},
$$
where  $\nx$ is the Marcinkiewicz space defined by the function \linebreak
 $\varphi_\varepsilon(t)=t\log_2^{-\varepsilon+1/2}(2/t)$.
\par
\medskip
 Proof}. Assume that the theorem is not true, that is,
for a certain  $\varepsilon>0$,
$$
\bar{T}_\theta :\,{\bar{\cal R}}(\de)\to{\nx}
\leqno{(33)}
$$
for any arrangement of the signs $\theta$.

Let
$$
x(s,t)=\ef\,\in{\,\de(\ab)},\,a=\pq\in{l_2}.
$$
By Lemma 3, for arbitrary natural $n$, the absolute values of the functions
$$
x_n(s,t)=\,\sum_{i,j=1}^na_{i,j}r_i(s)r_j(t)\,\,\,\mbox{and}\,\,\,y_n(t)=\,\sum
_ { i,j=1}^na_{i,j}r_i(t)r_{j+n}(t)$$ are equimeasurable.

Let us introduce the following arrangement of the signs
$\bar{\theta}=\{\bar{\theta}_{i,j}\},$
\linebreak $\bar{\theta}_{i,j}=\theta_{i,j-n}$,
if  $j>n$, and	$\bar{\theta}_{i,j}$ arbitrary, if  $j\le n$.
Then the absolute values of the functions
$$
T_{\theta}x_n(s,t)=\,\sum_{i,j=1}^n\theta_{i,j}a_{i,j}r_i(s)r_j(t)$$
and
$$
\bar{T}_{\bar{\theta}}y_n(t)=\,\sum_{i,j=1}^n\bar{\theta}_{i,j+n}a_{i,j}r_i(t)r_{j+n}(t)\,=\,\sum_{i,j=1}^n\theta_{i,j}a_{i,j}r_i(t)r_{j+n}(t)$$
are also equimeasurable.

Summarizing, from relations (33) we obtain for the arrangement
$\bar{\theta}$:
$$
\|T_{\theta}x_n\|_{\nx(\ab)}\,=\,\|\bar{T}_{\bar{\theta}}y_n\|_{\nx}\,\le C\|y_n\|_{\infty}\,=\,C\|x_n\|_{\infty}.$$
Since the norm in the spaces  $L_{\infty}$ and	$\nx$
is order semi-continuous
(see   \cite{9}), a passage to the limit as
$n\to{\infty}$ implies
$$
\|T_{\theta}x\|_{\nx(\ab)}\,\le\,C\|x\|_{\infty},$$
that is, for every arrangement of signs $\theta$,
$$
\wg:\,{\cal R}(\de)\,\to\,{\nx(\ab)},$$
which contradicts Theorem 7.
\par
\medskip
{\bf Corollary 3.} {\it If  $p\in{[1,2)}$,
then there is no s.s.  $X$ on  $I$ for which
$$
\biggl|\biggl|\sum_{i,j=1}^{\infty}b_{i,j}r_ir_j\biggr|\biggr|_X\,\asymp{\,\biggl(\sum_{i,j=1}^{\infty}|b_{i,j}|^p\biggr)^{1/p}}\,(b_{i,i}=0).$$}

The proof is similar to that of Corollary 2, with the only difference that
instead of Theorem 7 one has to apply Theorem 8.
\par
\medskip
{\bf Remark 4.} Propositions similar to Corollaries
 2 and	3 can be established for other s.s.
sequences, for example, for the spaces of Lorentz, Marcinkiewicz and Orlicz.
In the case of usual Rademacher system
the situation is completely different (see  \cite{3}, \cite{7}). For instance,
it was shown in   \cite{7} that for any space
of sequences
$E$, interpolation between   $l_1$ and	$l_2$,
one can find a functional s.s.	$X$ on	$[0,1]$
such that
$$
\biggl|\biggl|\fv\biggr|\biggr|_X\,\asymp{\,||(c_j)||_E}.
$$


\end{document}